\documentclass[journal]{IEEEtran}
\ifCLASSINFOpdf
 \usepackage[pdftex]{graphicx}

\graphicspath{{../images/},{../../SIMU/plots/keep/},{../Global_paper/images/}}
\else
\fi
\newif\ifIncludeFastProj 
\newif\iflongversion 

\longversiontrue

\IncludeFastProjtrue
\usepackage{amsmath}
\usepackage{frenchineq}
\usepackage{algorithm,algorithmicx}
\makeatletter
\renewcommand{\ALG@name}{Algo.}
\makeatother
\usepackage{algpseudocode}
\usepackage{cleveref}
\usepackage{bm}

\newcommand{\norm}[1]{\left\Vert #1\right\Vert}

\newcommand{\SC}{\text{S\hspace{-1pt}C}}
\newcommand{\rr}{\mathbb{R}} 
\newcommand{\argmin}[1]{\underset{#1}{\mathrm{arg}\hspace{-1.3pt}\min}}
\newcommand{\eqDef}{\vspace{-1pt}\overset{\Delta}{=}}

\usepackage[colorinlistoftodos,bordercolor=orange,backgroundcolor=orange!20,linecolor=orange,textsize=scriptsize]{todonotes}

\usepackage{multirow}
\usepackage{pbox}
\usepackage{amsfonts}
\usepackage{amssymb}
\usepackage{amsthm}
\usepackage{float}
\usepackage{adjustbox}
\usepackage{bm}
\usepackage{tikz}
\usetikzlibrary{arrows}

\usepackage[numbers,sort]{natbib}
\renewcommand{\citep}[1]{\citeauthor{#1} (\citeyear{#1})}
\newcommand{\citauty}[1]{\citeauthor{#1}, \citeyear{#1}}

\newtheorem{theorem}{Theorem}

\newtheorem{prop}{Proposition}

\newtheorem{rem}{Remark}
\newtheorem{assumption}{Assumption}

\newtheorem{definition}{Def.}

\Crefname{algorithm}{Algo}{Algos}
\Crefname{prop}{Prop.}{Props.}
\Crefname{equation}{Eq.}{Eqs.}
\Crefname{figure}{Fig.}{Figs.}
\Crefname{tabular}{Tab.}{Tabs.}
\Crefname{table}{Tab.}{Tabs.}
\Crefname{theorem}{Thm\@.}{Thms.}
\Crefname{definition}{Def.}{Defs.}
\Crefname{section}{Sec.}{Secs.}
\Crefname{assumption}{Assumption}{Assumptions.}

\newcommand{\G}{\mathcal{G}}
\newcommand{\BR}{\mathrm{B\hspace{-1.3pt}R}}
\newcommand{\NE}{\mathrm{N\hspace{-1.3pt}E}}

\newcommand{\dpart}[2]{\dfrac{\partial #1}{\partial #2}}
\newcommand{\x}{\ell}
\newcommand{\X}{\mathcal{L}}
\newcommand{\xx}{\bm{\ell}} 

\newcommand{\xag}{L} 
\newcommand{\xxag}{\bm{L}} 
\newcommand{\lbx}{\underline{\x}}

\newcommand{\ubx}{\overline{\x}}

\newcommand{\N}{\mathcal{N}}
\newcommand{\n}{n} 
\newcommand{\m}{m}
\newcommand{\T}{\mathcal{T}}
\renewcommand{\t}{t}
\newcommand{\ti}{_\t} 

\newcommand{\nt}{_{\n,\t}}
\newcommand{\mt}{_{\m,\t}}

\newcommand{\nf}{\mathrm{N\hspace{-1pt}F}}
\newcommand{\lnf}{L_{\nf}}
\newcommand{\lnft}{L_{\nf,t}}
\newcommand{\lnftp}{L_{\nf,t'}}

\newcommand{\Lnf}{\bm{L}_{\nf}}

\newcommand{\PoA}{\mathrm{PoA}}

\newcommand{\stronga}{a}
\newcommand{\lipCons}{M}

\newcommand{\ca}{\alpha}
\newcommand{\cb}{\beta}

\begin{document}
%
\title{
%
%
%
Analysis and Implementation of a Hourly Billing Mechanism for Demand Response Management
}
%
%
%

\author{Paulin Jacquot,~\IEEEmembership{Student Member,~IEEE,}
        Olivier Beaude, 
		St\'ephane Gaubert and Nadia Oudjane
\thanks{P. Jacquot, O. Beaude and N. Oudjane are with OSIRIS department in EDF Lab Saclay, France.}
\thanks{P. Jacquot and S. Gaubert are with Inria Saclay and CMAP, \'Ecole polytechnique, CNRS, Palaiseau, France.}
}

%
%

\markboth{
}%
{Shell \MakeLowercase{\textit{et al.}}: Bare Demo of IEEEtran.cls for IEEE Journals}
%



\maketitle


\begin{abstract}An important part 
of the Smart Grid literature on residential Demand Response  deals with  game-theoretic consumption models. Among those papers, the hourly billing model is of special interest as an intuitive and fair mechanism. We focus on this model and answer to several theoretical and practical questions. First, we prove the uniqueness of the consumption profile corresponding to the Nash equilibrium, and we analyze its efficiency by providing a bound on the Price of Anarchy. Next, we 
address the computational issue of the equilibrium profile by  providing two algorithms: the cycling best response dynamics and a projected gradient descent method, and by giving an upper bound on their convergence rate to the equilibrium. Last, we simulate this demand response framework in a stochastic environment where the parameters depend  on forecasts. We show numerically the relevance of an online demand response procedure, which reduces the impact of inaccurate forecasts.
\end{abstract}


%
\IEEEpeerreviewmaketitle

\section{Introduction}
%
%
%
%
\IEEEPARstart{D}{emand} Response (DR) is a technique to exploit electricity consumers flexibilities by giving them particular incentives, in order to achieve some services to the grid e.g. reducing production and distribution costs or increasing renewable energy insertion~\cite{siano2014demand}. In DR  programs, the aggregated energy demand is a key metric and an aggregator interacts with active consumers, willing to minimize their electricity bill or maximize their utility, to optimize this demand profile. In such a framework, energy can be viewed as an asset demanded by customers, and which has a cost that depends on total demand and the time of demand---some congestion effects arise on the most demanded time periods and as producers 
 have a limited quantity to offer. Game-theory offers a well adapted environment to model these congestion effects, and different game-theoretic framework have been proposed in the Smart Grid literature \cite{chai2014demand,chen2014autonomous,saad2012game}, following the seminal paper~\cite{mohsenian2010autonomous}. 
 Among these, the hourly billing mechanism, introduced in \cite{baharlouei2012tackling,baharlouei2013achieving} and also adopted in \cite{atzeni2013demand,chen2014autonomous} emerges as  a natural model---it has a structure of a \emph{routing congestion game} \cite{orda1993competitive}---enjoying important fairness properties, although its theoretical analysis is harder than the mechanism proposed in \cite{mohsenian2010autonomous}.  
  In  game-theoretic models, a major issue is to define an efficient procedure to compute and reach the consumption \emph{equilibrium} associated with the game. Several papers \cite{baharlouei2017convergence,atzeni2013demand,chen2014autonomous}  have investigated the complexity and algorithmic aspects associated to the notion of equilibrium. 
One can refer  to \cite{vardakas2015survey} for a survey of classical optimization methods and their applicability depending on the DR framework. In \cite{atzeni2013demand}, the authors consider  a distributed generation environment and use an algorithm based on an iterative proximal best response.  In \cite{chen2014autonomous}, a hourly billing mechanism is used and the authors propose a proximal-point algorithm to compute the equilibrium. 

 In this paper, 
 we investigate the theoretical properties and computational aspects of the hourly billing mechanism 
 and discuss its practical implementation. The main contributions are the following. First, we prove the uniqueness of the equilibrium (\Cref{th:uniqueness-NE}) under a convexity assumption. The uniqueness of the equilibrium profiles was proved for the daily billing mechanism proposed in \cite{mohsenian2010autonomous}. For the hourly billing mechanism, \cite[Prop\@. 1]{chen2014autonomous} gives the uniqueness for a particular class of price functions.  Our result applies to any convex and increasing price functions, and extends \cite[Thm.\ 1]{orda1993competitive} to a more general model where we consider bounds on the load at each time period. 
 Next, we give a bound on the resulting \emph{Price of Anarchy} (PoA)
  , which shows the efficiency of the equilibrium (\Cref{th:bound_PoA}). This result should be compared with \cite{mohsenian2010autonomous} where the model induces a PoA equal to one (optimality). Here, the PoA is numerically close to one but not one.
Then, we bound  the convergence rate of the best response algorithm in the case of affine prices (\Cref{th:cvg-BR-linear}). In that  case, convergence is known but to our knowledge, no bound on the rate has ever been given. The convergence has been conjectured more generally for any convex prices \cite{mertzios2008fast,brun2013convergence}. 
  We introduce a different algorithm based on a simultaneous projected gradient descent (\Cref{algo:SIRD}), and show its geometric convergence (\Cref{th:SIRD-convergence}) with a condition on the price functions only. To our knowledge, those results are also new. They can be compared with the algorithm 1 proposed in \cite{chen2014autonomous}. In particular, we allow a fix step-size and we do not need a proximal term. 
 Last, we introduce an online DR procedure with receding horizons, to take into account updated forecasts in a stochastic environment. We show numerically, based on  real consumption data, that this procedure can achieve significant savings compared to an offline procedure. 


This paper reassembles and extends the main results on the hourly billing model announced in our conference papers~\cite{Paulin2017ISGT,Paulin2017SGC}: 
 in \Cref{th:bound_PoA} here, we give the upper bound on the Price of Anarchy \cite[Thm.\ 2]{Paulin2017ISGT} and  in \Cref{prop:stable game HLP} we use the same property than \cite[Thm.\ 1]{Paulin2017ISGT}. We also use the potential property of the game noticed in \cite[Thm.\ 2]{Paulin2017SGC} and the \textit{Best Response} algorithm presented in \cite[Def.\ 3]{Paulin2017ISGT}. However there are several 
additional results in this paper: the theorem  of uniqueness of the equilibrium presented here (\Cref{th:uniqueness-NE}) is stronger than \cite[Thm.\ 1]{Paulin2017ISGT}. Also, the SIRD algorithm  (\Cref{algo:SIRD}) and the convergence theorems (\Cref{th:cvg-BR-linear,th:SIRD-convergence}) were not presented in \cite{Paulin2017ISGT,Paulin2017SGC}. We complete the simulation framework in \cite{Paulin2017ISGT,Paulin2017SGC} by considering  multiple forecasts on the nonflexible load and by introducing an online demand response procedure (\Cref{algo:onlineBRD}).

%

This paper is organized as follows: \Cref{sec:model} gives the mathematical model of the DR framework and the associated billing mechanism, under the form of a game.  In \Cref{sec:computation},  we introduce the notion of strong stability for general games and define two decentralized algorithms that enable to compute the equilibrium consumption profiles. We prove the convergence of those algorithms and give guarantees on their convergence rate. Finally, in \Cref{sec:simulations}, we define an \emph{online} DR procedure and simulate it with historical consumption data of consumers with electric vehicles as flexible consumptions. We compare the performance of this online DR scheme to the \emph{offline} version and other consumption scenarios.

\noindent\textsc{\textbf{Notation Convention}:} through this paper, bold font $\xx$ is used to denote a vector as opposed to a scalar $\x$. 

\section{Consumption Game with Hourly Billing}
\label{sec:model}
\newcommand{\hxx}{\hat{\xx}}

\subsection{District of flexible consumers}
We consider a network composed of a set $\N=\{1,\dots,N\}$ of residential consumers linked to a local \emph{aggregator}. Each household is equipped with a smart meter enabling two-way communication of information with the aggregator. 
We assume that each household electricity consumption can be divided into two parts: one which is \emph{nonflexible} (lights, cooking appliances, TVs) and one which is \emph{flexible} (Electric Vehicle charging, water heater, washing machine, etc).  Moreover, each smart meter is linked to an \emph{Electricity Consumption Scheduler} (ECS) that can automatically optimize and schedule the consumption profile of the flexible appliances of the user, given the constraints set by the consumer and physical constraints of each appliance.

\subsection{From individual to aggregated consumption profiles}

In the DR program, we determine a consumption profile for each consumer on a finite time horizon $T$. In this study, we take $\T$ as a discrete set of time periods $\T=\{ 1,\dots T\}$. In the simulations, $\T$ will correspond to one day, and each time period $\t$ to one hour.
We assume that each consumer has a Nonflexible (NF) consumption profile; we denote by $\Lnf \in \rr^T$, indexed by $\T$, the aggregated  nonflexible load profile on the set of consumers $\N$. On top her nonflexible consumption, each consumer $n$ has a flexible consumption profile that we denote by $ \xx_n = \left(\x_{n,1},\cdots,\x_{n,T}\right) \in {\rr}^T$. 
The aggregated flexible load profile on the set of consumers is obtained as:
\begin{equation}\label{eq:agg-flex-load}
\xxag=(\xag\ti)_{\t\in\T} \in \rr^T \hspace{2mm} \text{with} \hspace{2mm} \forall t \in \T, L_t \eqDef \sum_{n\in\N} \x\nt  \ . 
\end{equation}

\subsection{Aggregator cost: a function of aggregated consumption}
The aggregator is himself linked to electricity providers and we consider that he faces a per-unit (of energy) price function\footnote{This price function can represent different objectives of the aggregator as the minimization of the distance to a targeted aggregated consumption profile ${\xxag}^*$, the minimization of providing costs or the maximization of self-consumption of (local) renewable production.} $\xag\ti\mapsto c_{\t}(\xag\ti)$ associated with each time period $\t \in \T$ for the flexible electricity demand  $\xag\ti$ given in \eqref{eq:agg-flex-load}, so that the cost of providing $\xag\ti$ at $t$ is $\xag\ti\times c_{\t}(\xag\ti)$.  
Price functions can depend on the nonflexible consumption profile $\Lnf$ because they model a global cost for the aggregator. In particular, we will assume in our model that the aggregator faces for each time period providing costs $C_t(.)$ that depend on the total load $\lnft + \xag\ti$. The price function $c_t(.)$ for the flexible part of consumption at $t$ can be rewritten as:
\begin{equation} \label{eq:price-from-total}
c_t(\xag\ti) \eqDef \frac{1}{\xag\ti} \Big[ C_t(\lnft+\xag\ti) - C_t(\lnft) \Big] \ .
\end{equation}
\begin{rem}Following \eqref{eq:price-from-total}, the price function $c_t(.)$ depends implicitly on the value of the nonflexible load $\lnft$ (even if the dependancy is not explicit in the notations used hereafter).
\end{rem}
In our framework, we will consider the following assumptions on the price functions  $(c_t)_t$:
\begin{assumption}\label{ass:basic-properties-prices}
For each $t\in \T$, $c_t$ is twice differentiable, strictly increasing and convex.
\end{assumption}

\begin{assumption}\label{ass:strong-stability-condition}
For each $t\in \T$, $c_t$ is twice differentiable, convex and strictly increasing. Moreover, there exists $\stronga >0$ s.t. for any $t$ and admissible $\xx$:
\begin{align} \label{eq:strong-stability-assumption}
2c_\t'(\xag\ti)  \left( 1- \left(\textstyle\frac{c_\t''(\xag\ti)}{2c'_t(\xag\ti)}\right)^2 \norm{\xx\ti}_2^2 \right)\geq \stronga .
\end{align}
\end{assumption}

\begin{assumption}\label{ass:PricesLinear}
For each $t\in \T$, $c_t$ is affine, positive and increasing: $\forall t \in \T, \ c_t(\x)=\ca\ti+\cb\ti \x$ with $\ca\ti, \cb\ti \in (\rr^*_+)^2$ .
\end{assumption}
\begin{rem}
The latter three assumptions are more and more restrictive: \Cref{ass:PricesLinear} implies \Cref{ass:strong-stability-condition} with $\stronga=2\min_t \cb\ti$, and \Cref{ass:strong-stability-condition} obviously implies \Cref{ass:basic-properties-prices}.
\end{rem}
\begin{rem} For $\stronga=0$, inequality \eqref{eq:strong-stability-assumption} in \Cref{ass:strong-stability-condition} simplifies to the condition: $ \norm{ \xx\ti }_2^{-1} \geq \left| \frac{c''_\t(\xag\ti)}{2c_\t'(\xag\ti) }\right| \ $. For each $t$, $c_t''$ has to be small relatively to $c_t'$.\end{rem}

\Cref{ass:basic-properties-prices} is standard in the congestion games literature and corresponds to "type-B" functions in the seminal paper~\cite{orda1993competitive}. This assumption is also made in most of the papers dealing with game-theoretic DR models as \cite{baharlouei2013achieving}. Indeed, it is justified by the fact that marginal costs of producing and providing electricity are increasing. \Cref{ass:PricesLinear} is also a standard assumption made in \cite{baharlouei2014efficiency} because it enables fast computation of NE (see \Cref{sec:computation}), but it is very restrictive. \Cref{ass:strong-stability-condition} is not very explicit but is an in-between condition that comprises a larger set of functions than linear functions and for which our main results hold.

Whatever the assumption retained, the aggregator can influence the flexible consumption by sending incentives to consumers through a \emph{billing mechanism}, that is, by defining what each consumer $n$ will pay relatively to her flexible consumption profile $\xx_n$, and in turn which quantity $n$'s ECS will minimize.

\subsection{Consumer's Optimization Problems}
In this paper, following our studies in \cite{Paulin2017ISGT,Paulin2017SGC}, we will use an hourly proportional billing mechanism, where each consumer $n \in \N$ minimizes her bill:
\begin{equation}
\label{eq:bill_definition}
 b_n(\xx_n,\xx_{-n}) \eqDef \sum_{\t \in \T} \x\nt c_t(\xag\ti) \ ,
\end{equation}
where $\xx_{-n}\eqDef(\xx_m)_{m\neq n}$ denotes the consumption of all consumers but $n$.
This billing mechanism was shown to have interesting fairness properties and is also adequate when considering consumers' utility functions (representing, e.g., temporal preferences for flexible consumption) \cite{baharlouei2014efficiency,Paulin2017ISGT,Paulin2017SGC}. It corresponds to a routing ``atomic splittable'' congestion game framework \cite{orda1993competitive}, well studied in the game theory literature, where we add the bounding constraints \eqref{cons:minmax-power} below.

Through her ECS, each consumer will adjust her flexible consumption profile $\xx_n \in \rr^T$ to minimize her bill, which corresponds to the following optimization problem:
\begin{subequations}
\label{eq:user_problem}
\begin{align}
&\min_{\xx_n \in \mathbb{R}^T }  b_n(\xx_n,\xx_{-n}) \\
\label{cons:total-power} \text{s.t. } & \hspace{0.5cm}   \textstyle\sum_{\t\in \T} \ell\nt = E_{n} \ , \\ 
 \label{cons:minmax-power}& \hspace{0.5cm} \lbx\nt \leq \x\nt 
\leq \ubx\nt , \forall \t \in \T \ .
  \end{align}
\end{subequations}
 Constraint \eqref{cons:total-power} ensures that the total energy  given to $n$ satisfies her daily flexible energy demand  over $\T$, denoted by $E_n$, that we assume fixed and deterministic in our work\footnote{$E_n$  can be set by the consumer, induced by the physical parameters of her appliances (battery capacity), or computed by learning the consumer's habits.}. Constraint \eqref{cons:minmax-power} takes into account the physical power constraints and the personal scheduling constraints (supposed given by the user to her ECS). Note that taking $\lbx\nt=\ubx\nt=0$ forces $\x\nt=0$ so that constraint \eqref{cons:minmax-power} includes in particular unavailability during some time periods. We will denote by $\X_n$ the feasible set of user $n$, given by constraints (\ref{cons:total-power}-\ref{cons:minmax-power}), and $\X\eqDef\X_1 \times \dots \times \X_N$ the Cartesian product of the feasible sets. 
As $b_n$ depends both on $\xx_n$ and $\xx_{-n}$ through $\xxag\ti$, we get a $N$-person minimization game  that we write in the standard form \cite{fudenberg1991game} as $\G \eqDef \left( \N, \X, (b_n)_n \right) .$ \\
 
In game-theoretic models, a desirable stability property is when each player $n$ has no interest to deviate unilaterally from her current profile $\xx_n$. This corresponds to the notion of Nash Equilibrium ($\NE$), that is:
\begin{definition}[\citauty{nash1950equilibrium}]\textbf{ Nash Equilibrium (NE).}\\
$(\xx^\NE_n)$ is a $\NE$ of the minimization game:  $\G = \left( \mathcal{N}, \X, (b_n)_n \right)$ iff for any player $n \in \N$:
\begin{align*}
& b_n(\xx^\NE_n, \xx^\NE_{-n} ) \leq b_n(\xx_n, \xx^\NE_{-n}), \ \forall \xx_n \in \X_n  \\
\Leftrightarrow &  \ \xx^\NE_{n} \in \argmin{\xx_n \in \X_n} \ b_n(\xx_n, \xx^\NE_{-n}) 	\ . 
\end{align*}
\end{definition}
It is known that an NE may not exist or may not be unique, even in routing congestion games \cite{orda1993competitive}. In our framework however, both properties are ensured, as stated in \Cref{th:uniqueness-NE} below. This result extends the uniqueness theorem in \cite{orda1993competitive} in presence of the constraint \eqref{cons:minmax-power} on power bounds.
\begin{theorem}
\label{th:uniqueness-NE}
Under \Cref{ass:basic-properties-prices}, $\G$ has a unique NE.
\end{theorem}
\noindent
\iflongversion
\textit{Proof:} See Appendix~\ref{app:uniqueness-NE}.\vspace{0.2cm}
\fi

As said above, an NE is a very interesting situation in practice because of its stability: each player will only increase her bill if she changes her profile. However, an NE does not necessarily minimize the  \emph{social cost}: \begin{equation}
\SC (\xx) \eqDef\sum_n b_n(\xx) \ .
\end{equation}
 Note that, with the billing equation \eqref{eq:bill_definition}, $\SC(\xx)$ is equal to the total \emph{system cost}\footnote{In practice, the system costs can differ from the social cost of consumers, for instance if we consider that the aggregator makes a positive profit, or if we consider consumers utility functions as done in \cite{Paulin2017SGC}. } $\sum_t \xag\ti c_t(\xag\ti)$, so that this quantity should be minimized from the aggregator point of view. In general games, an NE can be sub-optimal in terms of \SC. To measure the inefficiency of Nash Equilibria, a standard quantity is the Price of Anarchy:
\begin{definition}[Koutsoupias \textit{et al}, \citeyear{koutsoupias1999worst}]\textbf{Price of Anarchy (PoA).}\\
Given a N-player game $\G=(\N,\X,(b_n)_n)$ and $\X_{\NE}$ its set of Nash equilibria, the PoA is defined as the following ratio:
\begin{equation*}
\PoA(\G)= \dfrac{\sup_{\xx \in \X_{\NE}} \SC(\xx) }{\inf_{\xx \in \X} \SC(\xx)} \ .
\end{equation*}
\end{definition} 
Note that, from the  definition, the $\PoA$ is always greater than 1. Furthermore, finding an upper bound on the $\PoA$ ensures that the social cost at any NE will be relatively close to the minimal social cost. In general, bounding the PoA is a hard theoretical question in general congestion games \cite{roughgarden2015intrinsic,johari2005efficiency}. In \cite{roughgarden2015local}, the authors give an upper bound if the price functions are polynomial with bounded degree and positive coefficients.  With degree one (affine prices, \Cref{ass:PricesLinear}) the bound is $\PoA\leq 1.5$, that is, 
 the NE profile can induce costs as much as 50\% higher than the optimal costs: implementing such a framework would not be worthwhile for the aggregator, as uncoordinated consumers will probably perform better (in our simulations, the uncoordinated profiles induce costs 16\% higher than the optimal costs, see  \Cref{tab:results}). However, the results in \cite{roughgarden2015local} are worst-case bounds and these bounds are only approached asymptotically\footnote{Meaning that there exists a sequence of games $(\G_\nu)_{\nu\geq0}$, with parameters depending on $\nu$, and affine price functions $c_t$ such that $\PoA(\G_\nu)\underset{\nu\rightarrow \infty}{\longrightarrow} 1.5$.}: in our simulations with affine prices, the PoA was always much lower than 1.5. One of the reasons is that in \cite{roughgarden2015local} the model does not consider the power constraints \eqref{cons:minmax-power}, and a PoA of $1.5$ might be reached in our case only if the constraints \eqref{cons:minmax-power} are coarse enough. To further explain the low PoA in our instances, we found the following theorem by precising the results of \cite{roughgarden2015local}: 
\newcommand{\ubxa}{\overline{\xag}}
\begin{theorem}\label{th:bound_PoA}
Under \Cref{ass:PricesLinear}, define for any $\t \in \T$, $\varphi_t=(1+\frac{\ca_\t}{\cb_\t \ubxa_\t})^2$, where $\ubxa_\t=\sum_n \ubx\nt$ and  $t_0 \eqDef \mathrm{arg}\min_{\t}\frac{\ca_\t}{\cb_\t \ubxa_\t}$. Assuming that, for all $t\in \T$:
\begin{equation}
\label{ass-prices-for-PoAbound}
\varphi_t \leq \varphi_{t_0} + 2 + \sqrt{1+\varphi_{t_0}} \ ,
\end{equation}
the following inequality holds:
\begin{align}
\label{eq:PoA-bound-with-phi}\PoA(\G) & \leq \frac{1}{2}(1+\sqrt{1+\varphi_{t_0}^{-1}} +\frac{1}{2}\varphi_{t_0}^{-\frac{1}{2}}) \ . 
\end{align}
%
\end{theorem}
\noindent
\iflongversion
\textit{Proof:} See Appendix~\ref{app:proof-PoAbound}.
\fi
\vspace{0.2cm}
\begin{rem}
Using the inequality $\forall x\geq0, \ \sqrt{1+x^2}\leq 1+x$, \eqref{eq:PoA-bound-with-phi} implies the following simplified---but coarser---bound:
\begin{equation}
\label{eq:PoA-bound-with-a-b}\PoA(\G) \leq 1+\textstyle\frac{3}{4}\sup_{\t\in\T}\left(1+\textstyle\frac{\ca\ti}{\cb\ti \ubxa_t}\right)^{-1}\ .
\end{equation}
\end{rem}
\noindent 
The assumption \eqref{ass-prices-for-PoAbound} in \Cref{th:bound_PoA} ensures that price functions $(c_t)$ cannot differ too much from one time period to another. This would be the case if, for instance, the price functions are uniform over $\T$. One can see that, according to \Cref{th:bound_PoA}, the PoA converges to one when $\ca\ti/(\cb\ti\ubxa_t) $ converges to infinity: the PoA can be arbitrarily close to one if we choose the coefficients $\ca\ti$ large enough. This result is indeed intuitive: when the prices are constant ($\cb\ti=0$), they do not depend on the aggregates $\xxag$ and there is no congestion effect; the optimal profile is obtained by each consumer choosing the time periods with lowest prices, independently of $\xx_{-n}$. Another interesting result is that the $\PoA$ also converges to one when the total load is low ($\ubxa_t\rightarrow 0$). Note that the left-hand-side of inequality \eqref{eq:PoA-bound-with-phi} is decreasing with $\varphi_0$ and is equal to $(\frac{1+\sqrt{2}}{2})^2\approx1.457$ for $\varphi_0=1$ so our result is always tighter than the bound given in \cite{roughgarden2015local}. However, in our simulations with linear prices, the PoA was still lower than the bound \eqref{eq:PoA-bound-with-phi}, even when assumption \eqref{ass-prices-for-PoAbound} does not hold: the inequality \eqref{eq:PoA-bound-with-phi} gives $\PoA\leq 1.271$ (average on the simulated days), while the PoA on mean values from \Cref{tab:results} is $1.017$. In this regards, getting a tighter bound or generalizing our proof to more general price functions could be the subject of future work.

\section{Fast Computation of the Nash Equilibrium}
\label{sec:computation}

The computation of NE is a central problem in game theory \cite{fabrikant2004complexity}. For a practical implementation of a DR program, we need to be able to compute the NE consumption profile in small time. In this section, we provide two algorithms for computing the NE and analyze their convergence in our specific setting.
\subsection{Two Decentralized Algorithms}

\label{subsec:computation}
\newcommand{\s}{\bm{s}}

Given a profile ${\xx}_{-n}$ of the others, consumer $n$ is expected to choose the profile $\xx_n$ corresponding to a minimizer of \eqref{eq:user_problem}, which is called its \emph{Best Response}\footnote{As player $n$ gives its \emph{best response} to the other players strategies.}. It is denoted by 
\begin{equation}\label{eq:best-resp-n}
\BR_n: \s_{n} \mapsto \argmin{\xx_n \in \X_n} \sum_t \x\nt c_t(s_{n,t}+\x\nt) \ ,
\vspace{-0.3cm}
\end{equation} 
which only depends on the sum of the load of the others $\s_{n} \eqDef \sum_{m\neq n } \xx_m \in \rr^T$ because of the ``aggregated'' structure \footnote{In a general setting, $\BR_n$ would be a function of $\xx_{-n}$.}: from \eqref{eq:bill_definition},  we see that $b_n$ only depends on $\xx_n$ and $\xxag_n=\s_n+\xx_n$.
A natural algorithm for computing an NE is to iterate best responses and update the strategies, cycling over the set of users until convergence. This procedure, described as \emph{Cycling Best-Response Dynamics} (CBRD) \cite{gilboa1991social} is described by  \Cref{algo:BRD}. 

\begin{algorithm}[H]
\begin{algorithmic}[1] 
\Require ${\xx}^{(0)}, \ k_{\rm{max}}, \ \varepsilon_{\rm{stop}}$
\State $k \leftarrow 0$, $ \  \varepsilon^{(0)} \leftarrow\varepsilon_{\rm{stop}}$
\While{ $ \varepsilon^{(k)} \geq \varepsilon_{\rm{stop}} $ $\And$ $k\leq k_{\rm{max}}$}
\For{$n=1$ to $N$} \label{line:algo-BRD-for}
\State $\s_{n}^{(k)}= \sum_{m<n} \xx_m^{(k+1)} + \sum_{m>n} \xx_{m}^{(k)}$ \label{line:algBR-update-s}\; 
\State $\xx_n^{(k+1)} \leftarrow \BR_n(\s_{n}^{(k)}) $ \hspace{3mm} (using \eqref{eq:best-resp-n})\label{line:algo-BRD-BR}\;
\EndFor
\State $ \varepsilon^{(k)}= \norm{ {\xx}^{(k+1)}-{\xx}^{(k)} }$ \;
\State $  k \leftarrow k+1 $ \;
\EndWhile
\end{algorithmic}
\caption{{Cycling Best Response Dynamics (CBRD) } }
\label{algo:BRD}
\end{algorithm}
\vspace{-0.5cm}
\begin{rem} In a general setting, $\BR_n(\s_{n})$ can be multivalued. In that case, we can still use \Cref{algo:BRD} by arbitrarily choosing any element of $\BR_n(\s_{n})$ at \Cref{line:algo-BRD-BR}. 
\end{rem}
Notice that the \textbf{for} loop in  \Cref{algo:BRD} (\Cref{line:algo-BRD-for}) implements sequential updates and cycles over the set of players in the arbitrary order\footnote{Choosing a "good" order of the $\BR$ in the \textbf{for} loop might accelerate the convergence of the algorithm. This could be the subject of future work.} $1,2,\dots,N$. One could think of a \emph{simultaneous} version of \Cref{algo:BRD} (without \Cref{line:algBR-update-s} and where \Cref{line:algo-BRD-BR} is executed by all players in parallel). However, we observed that doing so can prevent the convergence of \Cref{algo:BRD}.

Another natural algorithm to compute the equilibrium is to emulate the projected gradient descent, well-known in convex optimization \cite{calamai1987projected}, by considering the gradient of each objective function of the players, as described in \Cref{algo:SIRD}. 

\begin{algorithm}[H]
\begin{algorithmic}[1]
\Require ${\xx}^{(0)}, \ k_{\rm{max}}, \ \varepsilon_{\rm{stop}}, \ \gamma $
\State $k \leftarrow 0$, $ \  \varepsilon^{(0)} \leftarrow\varepsilon_{\rm{stop}}$
\While{ $ \varepsilon^{(k)}\geq  \varepsilon_{\rm{stop}} $ $\And$ $k\leq k_{\rm{max}}$}
\For{$n=1$ to $N$} \label{alg-line:for-loop}
\State $\xx_n^{(k+1)} \leftarrow \Pi_{\X_n}\left( 	\xx_n^{(k)} - \gamma \nabla_n b_n(\xx_n^{(k)},\xx_{-n}^{(k)}) \right)$ \label{line-alg:SIRD-proj}
\EndFor
\State $ \varepsilon^{(k)}= \norm{ {\xx}^{(k+1)}-{\xx}^{(k)} }$ \;
\State $  k \leftarrow k+1 $ \;
\EndWhile
\end{algorithmic}
\caption{{Simultaneous Improving Response Dynamics (SIRD) } }
\label{algo:SIRD}
\end{algorithm}
At \Cref{line-alg:SIRD-proj} of \Cref{algo:SIRD}, $\Pi_{\X_n}$ denotes the projection on the feasibility set of consumer $n$ $\X_n$.  The chosen denomination \emph{improving response} recall that, at each iteration of \Cref{algo:SIRD}, player $n$ \emph{improves} her profile $\xx_n$ by performing a projected gradient step (\Cref{line-alg:SIRD-proj}), but in general does not choose the best improvement as in \Cref{algo:BRD}.
Note that from \Cref{algo:BRD} to \Cref{algo:SIRD}, only the instructions within the \textbf{for} loop are changed: the update of $\s_n$ and computation of $\BR_n$ (\Cref{line:algBR-update-s,line:algo-BRD-BR} of \Cref{algo:BRD}) are replaced with the gradient step (\Cref{line-alg:SIRD-proj} of \Cref{algo:SIRD}).

\begin{rem}
Both \Cref{algo:BRD} and \Cref{algo:SIRD} can be implemented in a ``decentralized'' procedure: the instructions within the \textbf{for} loop (Line 4-5 for \Cref{algo:BRD} and \Cref{line-alg:SIRD-proj} for \Cref{algo:SIRD}) can be performed locally by each consumer's ECS. In this way, the privacy of consumers is respected as they do not have to send information about their constraints (\ref{cons:total-power}-\ref{cons:minmax-power}) to the aggregator. On the other hand, they only need to receive information on the aggregated load $\s_{n}^{(k)}$ and can hardly deduce the individual consumption of the other consumers.
\end{rem}
%

The computational complexity of one iteration (within the \textbf{for} loop) of \Cref{algo:SIRD} is equivalent to the complexity of the projection $\Pi_{\X_n}$,
\ifIncludeFastProj{which can be computed with the Quadratic Program (QP) $\Pi_{\X_n}(\xx'_\n)=\text{argmin}{}_{\xx_\n \in \X_n}\norm{\xx'_\n-\xx_n }_2^2$ so it would be of the same order of complexity (see \cite[Lecture 4]{ben2001lectures}) as one iteration (within the \textbf{for} loop) of algorithm CBRD with affine prices (\Cref{ass:PricesLinear}). With the specific structure of the feasible set $\X_n\subset\rr^T$, a QP can be solved very efficiently in $\mathcal{O}(T)$ \cite{brucker1984n}, so that one iteration of \Cref{algo:SIRD}---as well as one $\BR$ with linear prices in \Cref{algo:BRD}---will be very fast. Moreover, as we do not update sequentially the load of the others $\xx_{-n}$ in \Cref{algo:SIRD}, the projected gradient step within the \textbf{for} loop can be computed simultaneously and can be parallelized.
}
\fi
\subsection{Game Stability and Convergence of \Cref{algo:BRD,algo:SIRD} }
\label{sec:properties}

To study the convergence to the unique NE of the two algorithms proposed in \Cref{subsec:computation}, we use the notion of \emph{stability}, and prove that the energy consumption game $\G$ defined above is strongly stable under \Cref{ass:strong-stability-condition}.  The notion of stability was introduced in \cite{hofbauer2009stable}  in order to study different game dynamics in continuous time and their convergence to NE. We extend this property to a ``strong'' version (symmetrically to the concept of strong monotonicity for operators):

\newcommand{\xxp}{\xx'}
\begin{definition}[\citauty{hofbauer2009stable}]{\textbf{Stable Game.\\} }
\label{def:stable game}
\label{def:strongly-stable-game}
A minimization game $\G= \left( \N, \X, (b_n)_n \right) $ is stable iff 
\begin{equation} \label{eq:stability}
\forall \xx,\xxp \in \X, \ (\xxp-\xx)^\text{T} . \left(F(\xxp)-F(\xx) \right) \geq 0 \ ,
\end{equation}
with $F(\xx)\eqDef \left( \nabla_n b_n(\xx) \right)_{n\in\mathcal{N}} .$ 

\noindent Moreover, $\G$ is $\stronga$-strongly stable, with a constant $\stronga>0$, iff:
\begin{equation} \label{eq:strong-stability}
\forall \xx,\xxp \in \X, \ (\xxp-\xx)^\text{T}  . \left(F(\xxp)-F(\xx) \right) \geq {\stronga} \norm{ \xx-\xxp}^2 \ .
\end{equation}

\end{definition}

\begin{rem} The condition of stability in \eqref{eq:stability} is equivalent to the condition of strict diagonal convexity in \cite{rosen1965existence}, which implies uniqueness of $\NE$ \cite[Thm.2]{rosen1965existence}.
\end{rem}

%
\Cref{def:stable game} gives an abstract condition on an operator that depends on the objective functions of the players. In our case, players objectives $(b_n)$ depend linearly on price functions $(c_t)_t$ through
\eqref{eq:bill_definition}, so it is interesting to translate the  condition of \Cref{def:stable game} directly on the price functions, as stated in \Cref{prop:stable game HLP}.
\begin{prop}
\label{prop:stable game HLP}
Let $\stronga>0$ such that \Cref{ass:strong-stability-condition} holds. Then, the game $\G$ is $\stronga$-strongly stable.
\end{prop}
\noindent
\iflongversion
\textit{Proof:} See Appendix~\Cref{app:strong-stability}.
\fi

\vspace{0.5cm}
On the basis of this stability property, we first consider the question of the convergence of \Cref{algo:BRD}. In general games, CBRD might not converge \cite{hofbauer2016discretized} or might take an exponential time to converge \cite{ackermann2008convergence}. 
In atomic splittable congestion games on a parallel network, as in our case, the convergence and the speed of \Cref{algo:BRD} has been studied previously  in \cite{mertzios2008fast} and \cite{brun2013convergence}, where the authors show by different methods that there is a geometric convergence in the case of $N=2$ players and convex and strictly increasing price functions (\Cref{ass:basic-properties-prices}). However, to the best of our knowledge, the convergence in this setting and for more players $N>2$ is still an open question.

\begin{figure}[ht]
\vspace{-0.5cm}
\centering
\includegraphics[width=1.0\columnwidth]{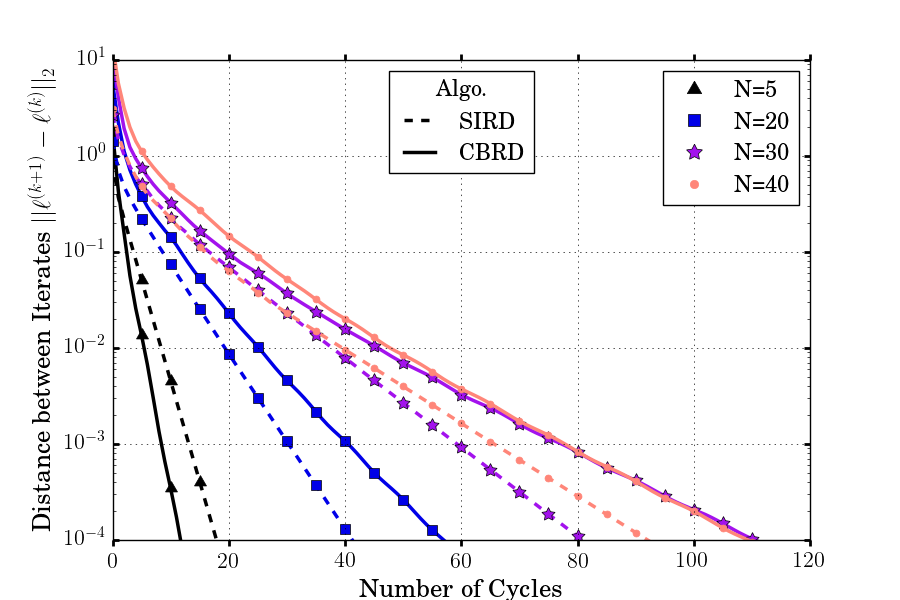}
\caption{Convergence of SIRD and CBRD with uniform affine price functions and $T=10$. \textit{When the number of players $N$ increases, the convergence rate of both algorithms decreases, but SIRD becomes faster than CBRD.}}
\label{fig:convergenceSIRDvsBR}
\vspace{-0.3cm}
\end{figure}

 In our case, simulations show a geometric convergence rate for any instance of $\G$ satisfying \Cref{ass:PricesLinear} and for any $N\in \mathbb{N}$, as illustrated in \Cref{fig:convergenceSIRDvsBR}. In \cite{brun2013convergence}, it is conjectured that this geometric convergence may also holds under \Cref{ass:basic-properties-prices}. Restricting ourselves to affine price functions, we notice that game $\G$ is a potential game \cite{Paulin2017SGC,monderer1996potential} and we get the following guarantee on the rate of convergence of \Cref{algo:BRD}:
 \begin{theorem}\label{th:cvg-BR-linear}
Under \Cref{ass:PricesLinear}, the sequence of iterates of Algorithm CBRD $\left(\xx^{(k)} \right)_{k\geq 0 }$ converges to the unique NE $\xx^\NE$ of $\G$. Moreover, the convergence rate satisfies:
\begin{equation*}
\norm{ \xx^{(k)}-\xx^\NE}_2 \leq C\frac{ \sqrt{\lipCons} N}{\sqrt{\stronga}} \times\frac{1}{\sqrt{k}} \ ,
\end{equation*}
where $C$ depends on $\x^{(0)}$ and the billing functions, \\ $\lipCons=2\max_t \cb\ti$ and $\stronga=2 \min_t \cb\ti  $.\end{theorem}
\noindent
\iflongversion
\textit{Proof:}  See Appendix~\ref{app:convergence-CBRD}.
\fi
The result is implied by convergence of alternating block coordinate minimization method \cite{hong2017iteration}.

The proof of \Cref{th:cvg-BR-linear} uses the fact that $M=\max_n M_n $ where $M_n$ is a Lipschitz constant of $\nabla_n b_n$, and $\stronga$ is a strong convexity (and $a$-strong stability) constant. To the best of our knowledge, the question to know if \Cref{th:cvg-BR-linear} holds for general price functions is open; it can be an avenue for future research.

It is easier to get a strong guarantee on the convergence rate of \Cref{algo:SIRD} for general price functions, as stated in \Cref{th:SIRD-convergence}:
\begin{theorem}\label{th:SIRD-convergence}
Denote by $\lipCons_n$ a Lipschitz constant of $\nabla_n b_n$ and  $\lipCons\eqDef \max_n \lipCons_n$.
 Under \Cref{ass:strong-stability-condition} ($\stronga$- strong stability), for $\gamma\eqDef \stronga/(N \lipCons^2)$, SIRD converges. Moreover, we have:
$$\norm{{\xx}^{\NE}-{\xx}^{(k)}}_2< \eta^k \norm{{\xx}^{\NE}-{\xx}^{(0)}}_2 \ ,$$
where $\eta=1-\frac{\stronga^2}{N\lipCons^2} \,$.
\end{theorem}

\noindent
\iflongversion
\textit{Proof:} See Appendix~\Cref{app:convergence-SIRD}.
\fi

Note that, under \Cref{ass:PricesLinear}, as stated in \Cref{th:cvg-BR-linear} we have $\lipCons=2\max_t \cb\ti$ and $\stronga=2 \min_t \cb\ti$, which gives the explicit contraction ratio $\eta= 1-\frac{\max_t \cb_t}{N \min_t \cb_t}$.

In practice, in spite of the weaker convergence result for CBRD, the convergence seems to also happen at a geometric rate, with a better ratio than the one found for Algorithm SIRD when the number of players $N$ is small, as illustrated in \Cref{fig:convergenceSIRDvsBR}. The convergence speed of both algorithm decreases with the number of users, as illustrated by the geometric coefficient $\eta$ in \Cref{th:SIRD-convergence}. However, SIRD becomes faster than CBRD when the number of players is large enough ($N \geq 20$).
\renewcommand{\figurename}{Fig.}

\section{Simulation of Online Demand Response }
\label{sec:simulations}

\newcommand{\lnfFor}{\hat{\xag}_{\nf}}
\newcommand{\lnfFort}{\hat{\xag}_{\nf,\t}}
\newcommand{\lnfForT}{\hat{\xag}_{\nf,T}}
\newcommand{\lnfFortp}{\hat{\xag}_{\nf,\t'}}

\newcommand{\LnfFor}{\hat{\xxag}_{\nf}}
\newcommand{\Ltot}{\bm{L}_{\text{tot}}}
\newcommand{\ltot}{L_{\text{tot}}}

\newcommand{\tend}{T }
\newcommand{\tst}{1}
In this section, we simulate consumption under the DR framework described above. We propose a practical procedure to implement the DR framework: at each hour, the equilibrium profiles for the flexible consumption is re-computed for the hours ahead to the end of the optimization horizon $T$, using \Cref{algo:BRD} or \Cref{algo:SIRD}. Next, we detail this simulation framework. 

\subsection{Online Demand Response Procedure}
\newcommand{\tmvi}{^{(t)}}
\label{subsec:DR-procedure}
The initial time horizon $\T$ that we consider for the planning via DR starts each day at noon ($t=1$) and stops at noon the day after ($t=T$). We describe an ``online'' procedure, computing the DR equilibrium flexible consumption profiles on time horizon $\{1,\dots T \}$ for each day. As the price functions $c_t$ depend on the nonflexible load through \eqref{eq:price-from-total}, and as the accuracy of forecast of this load improves when approaching from real-time, we re-compute the equilibrium using updated forecasts  at each time period.

\begin{algorithm}[H]
\begin{algorithmic}[1]
\State Start at $t=1$
\While{ $ t \leq T $}
\State Set new horizon $ \T\tmvi =\{t,t+1,\dots, T \}$ \;
\State Get  $\Lnf$ forecast on $\T\tmvi$:  $\LnfFor\tmvi\eqDef  \big({\lnfFor}\tmvi{_{,s}} \big)_{t\leq s\leq T }$
\label{line:get-forecast}\;
\State Re-compute prices $c_t(.)$ for $t\in \T\tmvi $ from \eqref{eq:price-from-total} \label{line-alg:onlineBRD-compute-prices}\;
\State Run Algo. SIRD or BRD to compute NE ${\xx}\tmvi$ on $\T\tmvi$\label{line:NEcompute}\;
\For{each user $n\in\N$ }
\State Realize computed profile on time $t$,  $\x\tmvi\nt$ \label{line:realization} \;
\State Update energy demand $E_n \leftarrow E_n-\x\tmvi\nt$ \label{line:energyUpdate}\;
\EndFor
\State Wait for $t+1$ \;
\EndWhile
\end{algorithmic}
\caption{{Online Demand Response Procedure} }
\label{algo:onlineBRD}
\end{algorithm}
\vspace{-0.4cm}

%
\begin{rem}
In practice, the NE profile $\xx\tmvi$ has to be computed before period $t$ to begin consumption at time $t$ (\Cref{line:realization}). If $\tau$ is an upper bound on the computation time of the NE profile (\Cref{line:NEcompute}), then, as we want to use the latest available forecast, Lines 3-5 would be run just before $t-\tau$, \Cref{line:NEcompute} is run in the interval $[t-\tau,t]$ and \Cref{line:realization} is executed through $[t,t+1]$.
\end{rem}
Observe that in \Cref{algo:onlineBRD}, considering updated forecast at \Cref{line:get-forecast} leads to updated price functions $(c_t)_t$ (\Cref{line-alg:onlineBRD-compute-prices}), according to equation \eqref{eq:price-from-total}. In turn, the updated price functions modify the objective function of user $n$, $b_n$, used in \Cref{line:NEcompute}.

The difference of \Cref{algo:onlineBRD} with an ``offline'' version is that we recompute the equilibrium consumption (\Cref{line:NEcompute}) at each time for all the time periods ahead. In an offline DR, we would compute the equilibrium consumption for all the horizon $\T=\{1,\dots,T \}$ \emph{only once}, just before $t=1$.

Proceeding with this ``online'' version has two main advantages. First, it enables to rely on updated forecast with new information acquired on the nonflexible load $\Lnf$ (\Cref{line:get-forecast}) up to time $t-\tau$. Second, it also enables to cope with local issues as disconnection of an user or a communication bug: in that case, we do not follow lines \ref{line:realization} and \ref{line:energyUpdate} for the involved user, and this user will have the same energy demand for the next round at $t+1$.

\begin{theorem}
\label{th:consistency-DR-proc}
The online demand response procedure of \Cref{algo:onlineBRD} is consistent: that is, if forecasts are perfect (i.e. $\forall t \in \T, \forall t' \in \T\tmvi, \lnfFor\tmvi{}_{,t'}= \lnf{}_{,t'}$), then for any $t_2>t_1$, the NE profile $\xx^{(t_1)}$ computed at $t_1$ with forecast $\LnfFor^{(t_1)}$ is equal on $\{t_2,\dots \tend\}$ to the NE profile $\xx^{(t_2)}$ computed at $t_2$ with forecast $\LnfFor^{(t_2)}$.
\end{theorem}
\noindent
\iflongversion
\textit{Proof:} See Appendix~\ref{app:proof-DR-procedure}.
\fi

\Cref{th:consistency-DR-proc} states a \emph{dynamic programming principle} adapted to our game-theoretic framework. It ensures that, following \Cref{algo:onlineBRD}, the final realized profile $\xx$ will correspond to the NE under perfect forecasts. 
To quantify the value of this online procedure in the more realistic case of imperfect forecasts, we simulate it on the set of consumers and parameters taken from real data, defined below.

\subsection{Consumers}
\newcommand{\D}{\mathcal{D}}
\newcommand{\nd}{_{n,d}}
\newcommand{\ndt}{_{n,d,t}}
\newcommand{\dt}{_{d,t}}
We consider a set of $N=30$ users owning an electric vehicle (EV) from the database of Texan residential consumers \textit{PecanStreet Inc.} \cite{PecanStreet}. We consider that the charging of the EV is the only flexible appliance of the consumers managed through the DR program, while the remaining of the user's consumption is nonflexible and is taken as in the data. 
We denote by $\D\eqDef \{ 16/01/01, \dots ,16/01/31\}$ the set of the 31 days of January 2016 for which we simulate the DR program and we index a parameter by $d\in \D$ when it is specific to day $d$. For constraints (\ref{cons:total-power}-\ref{cons:minmax-power}), we take, for each day $d\in\D$, the total flexible demand of user $n$, $E\nd$ as the total observed consumption for the EV of $n$ on the time set $\T=\{1,\dots,T\}$ (twenty-four hours of $d$ from 12\textsc{pm} to 11\textsc{am} on $d+1$). The power lower bound is always taken to zero $\underline{\x}\ndt =0$. For the power upper bound $\ubx\ndt$, we consider two cases: if a positive power was given at $d,t$ in the data, we   
take $\overline{\x}\ndt$ as the maximum power given to $n$'s EV in the data in the set $\D$. If the power given to the EV is 0 at $d,t$ in the data, we take $\overline{\x}\ndt=0$ (\textit{i.e.} we consider that the EV of $n$ was not available to charge on period $d,t$). 

\subsection{Price Functions}
\newcommand{\totalDemand}{D}
Following \cite{Paulin2017SGC}, we consider that the aggregator has a providing cost for the global demand at time $t$,  $\totalDemand\ti \eqDef(\lnft+\xag\ti)$, that does not depend on the time, and given (in \$) by 
$ C(\totalDemand\ti)
\eqDef 0.711-0.0417\totalDemand\ti+0.00295\totalDemand\ti^2  $
\newcommand{\avg}{\operatorname{avg}}
where the coefficients replicate the cost function of a real residential electricity provider. For this, we computed the average, minimum and maximum values of $\lnft$ over all the hours of the 31 days of January 2016 on our set of 30 consumers and interpolate the three values (avg $\lnft $, $ \min \lnft$, $ \max \lnft$) to three respective prices proposed by the Texan distributor \textit{Coserv} \cite{Coserv}  so that $\tilde{c}(\avg \lnft)= 0.080\$ $/kWh (price for ``base'' contracts), $\tilde{c}(\min \lnft)= 0.055\$ $/kWh (price for Off-Peak hours  in Time-of-Use contracts) and $\tilde{c}(\max \lnft)= 0.14\$ $/kWh (price for Peak hours).
Following \eqref{eq:price-from-total}, the price for the flexible load is given by: 
$  c_t(\xag\ti)
= (-4.17+0.590\lnft) + 0.295\xag\ti $, so that \Cref{ass:PricesLinear} holds.


\newcommand{\hl}{\hat{L}}
\newcommand{\dr}{\mathrm{d}}

\subsection{Forecast of the Nonflexible Load}

\begin{figure}[!ht]
\vspace{-0.5cm}
\centering
\includegraphics[width=1.0\columnwidth]{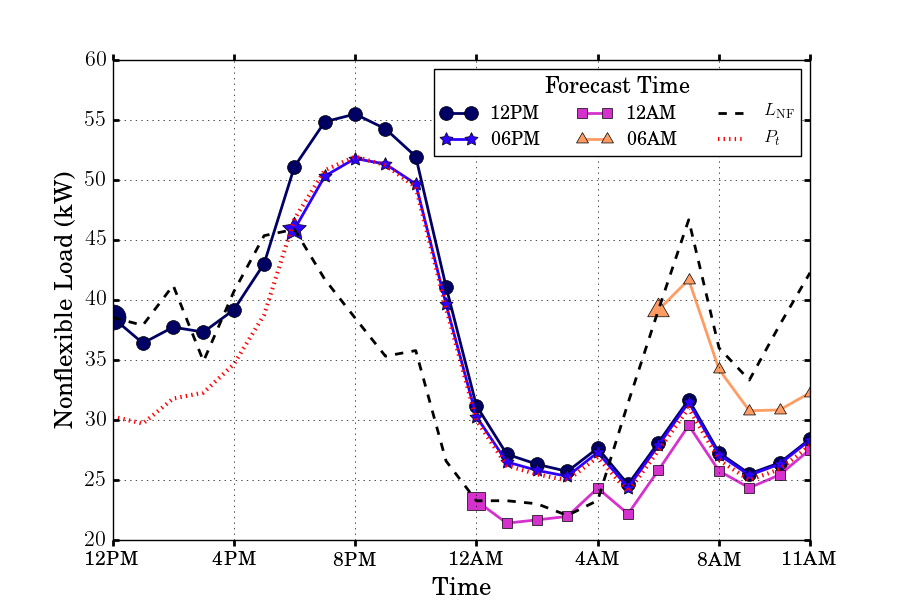}
\caption{Forecasts of the nonflexible load $\LnfFor\tmvi$  evolving in time. \newline \textit{We assume a perfect forecast at time} $t$ \textit{for }$t$. \textit{Forecasting performance increases when approaching real time.}}
\label{fig:forecasts}
\vspace{-0.2cm}
\end{figure}

\newcommand{\drift}{m}
Here, as we assume that the prices depend on the nonflexible load, at each time $t$ the aggregator has to compute a forecast  $\LnfFor\tmvi\eqDef  \big(\lnfFort\tmvi , \dots \lnfForT\tmvi\big)$ 
to be able to compute the equilibrium consumption for time periods $\{t,\dots, \tend\} $.
%
 To simulate the forecasts, we assume that the forecast made at time $t$ for period $t'\geq t$, $\lnfFortp\tmvi$ has no bias, that is $\mathbb{E}[\lnfFortp\tmvi  | \sigma(\mathcal{F}_t) ]=\lnftp$ (where $\mathcal{F}_t$ is the natural filtration over $(\lnft)_t$), and that we have perfect information at time $t$, that is: $\lnfFort\tmvi=\lnft \, .$ Considering that $\lnft=P_te^{X_t}$ where $X_t$ follows an Ornstein-Uhlenbeck  \cite{Ornstein1930} process with mean reversing coefficient $\drift$ and volatility $\sigma$, and $P_t$ a seasonality factor that depends on the hour of the week ($1^{\text{st}}$ hour to $168^{\text{th}}$ hour), we get for any $\t\leq t'$:
\begin{align*}
  \lnfFortp\tmvi= & P_{t'} \left( \frac{\lnft}{P_t} \right)^{e^{-\drift({t'}-t)}} \hspace{-5pt}\text{exp}\left( \frac{\sigma^2}{4\drift}(1-e^{-2\drift({t'}-t)}) \right) .  
\end{align*}
Using a least-squares regression on the observed data from years 2014 and 2015, we compute $\drift\simeq 0.198$ h$^{-1}$ and $\sigma\simeq 0.117$ h$^{-1/2}$. 
 An example of the simulated forecasts made at four different time periods is given in \Cref{fig:forecasts}.
 \vspace{-0.5cm}

\subsection{Gains with the Online DR Procedure} \label{subsec:results}

We run the DR Procedure as described in \Cref{subsec:DR-procedure} on each day of January 2016 to get the computed flexible consumption profile of each user $\xx_n$, and compute the associated social cost on the DR horizon $\{\tst,\dots, \tend \}$ for each day in $\D$.
We compare the total social costs on all simulated days $\D$, obtained via the DR online procedure and the associated social costs with the four other consumption scenarios below:

1) \emph{uncoordinated} case: no DR is implemented to control  or incentivize consumers flexibility; the consumption profiles are taken as the observed value in the data;

2) \emph{offline DR}: the equilibrium is computed only once at $t=\tst$ and for the whole time horizon $\{\tst, \dots \tend\}$ considering the first forecast $\LnfFor^{(\tst)}$ available at $t=\tst$,

3) \emph{perfect forecast DR}: offline DR, where we take $\LnfFor^{(\tst)}=\Lnf$. With \Cref{th:consistency-DR-proc}, it is useless to recompute the profiles at each time period,

4) \emph{optimal} scenario: a centralized entity (with perfect forecasts) computes the flexible consumption profile $\xx$  that minimizes the system costs $\sum_t \xag\ti c_t(\xag\ti)$.

The NE is computed by implementing \Cref{algo:BRD} under Python 3.5, where each quadratic program (QP) minimization (corresponding to one \emph{Best Response}, \Cref{line:algo-BRD-BR}) is solved with the solver Cplex 12.6 through Pyomo API, run on a single thread on an \textit{Intel i7}@2.6GHz. As for the stopping criterion of \Cref{algo:BRD}, we take $\varepsilon_\text{stop}=10^{-3}$. 
Under this configuration it takes on average around 80 seconds to compute each NE ($\xx_n^{\NE} \in \rr^{24}, \ n \in \{1\dots 30\}$).
The optimization problem to compute the \emph{optimal} consumption profile satisfying all users constraints (\ref{cons:total-power}-\ref{cons:minmax-power}), for each simulated day (from 12\textsc{pm} to 11\textsc{am}) in $\D$, is also a convex QP that is solved easily with the solver Cplex 12.6 in 0.31seconds on average.

\begin{table}[H]
\centering
\begin{tabular}{|c|c|c|c|}
\hline
\textbf{Cons. Scenario} & \textbf{Social Cost} & \textbf{Avg. Price} & \textbf{Gain} \\\hline\hline
Uncoordinated & \$ 1257.2 & 0.200 \$/kWh & ---\\\hline
Offline DR &  \$ 1231.6 & 0.195 \$/kWh &2.036\% \\\hline		
\textbf{Online DR} & \textbf{\$ 1131.1} & 0.180 \$/kWh &\textbf{10.03\%} \\\hline 
Perfect forecast DR & \$ 1075.2 & 0.171 \$/kWh & 14.47\% \\\hline 			
Optimal scenario & \$ 1056.8 & 0.169 \$/kWh & 15.94\%  \\\hline
\end{tabular}
\vspace{2pt}
\caption{Social Costs, average prices and relative gain to the uncoordinated consumption scenario on January 2016 .}
\label{tab:results}
\end{table}
\vspace{-0.5cm}

\begin{figure}[!ht]
\vspace{-0.5cm}
\centering
\includegraphics[width=1\columnwidth]{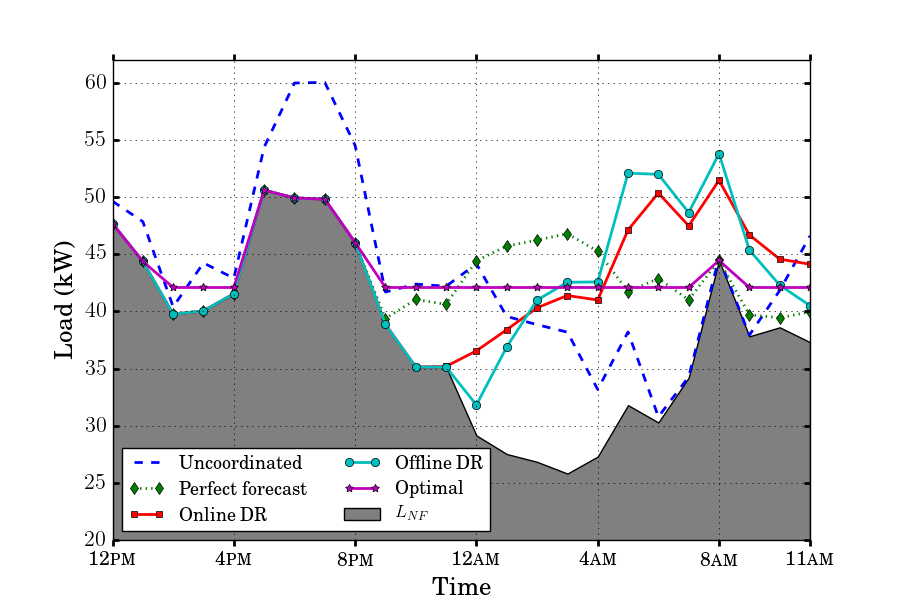}
\caption{Consumption Profiles on a typical day, with the different scenarios listed in \Cref{subsec:results}. \textit{The optimal profile flattens the consumption. The online DR procedure of \Cref{algo:onlineBRD} gets closer to the Perfect forecast (offline) DR profile.}}
\end{figure}

\Cref{tab:results} summarizes numerical results, it gives the total costs on the 31 days of January 2016 and compares the gains of the different flexible consumption scenarios relatively to the \emph{uncoordinated} one. We can see on this table that the online DR procedure achieves significant savings compared to the offline version for which the results are really low in average: using our offline DR decreases the system costs by $2\%$ relatively to the uncoordinated profile, that is, when consumers behave without any incentives.
 Implementing this offline DR program might not be worthy as it still involves a sophisticated communication and automation structure and it adds more constraints for consumers. This low performance is directly linked to our simple and naive model for the nonflexible load forecasts, which results in inaccurate forecasts for the last hours, as seen in \Cref{fig:forecasts}. If more advanced forecasting methods (see \cite{alfares2002electric}) can improve this accuracy, we cannot get rid of the high variance due to the low scale of the small set of consumers (30 in our example, and several hundreds for an aggregator for a typical low-voltage station).
 The online DR procedure seems to bring a solution to this issue: even with our simple forecast model, we achieve more than $10\%$ of savings, with a gap of only $5\%$  from the scenario with perfect forecasts. Therefore, our results show that implementing the given online DR procedure, even without very accurate forecasts, will be worthwhile for the aggregator.

\section{Conclusion}
In this paper, we developed a game-theoretic model for a residential demand response program, and we have addressed several issues about its implementation. We gave several new theoretical results about the uniqueness and existence of a Nash equilibrium consumption profile, and we have shown that the two algorithms CBRD and SIRD provide approximations of the NE at an arbitrary accuracy in finite time. We have introduced and simulated an online procedure that recomputes the NE profiles at each time period to take into account new information. We have shown numerically that this online procedure achieves a small price of anarchy when the parameters are fixed but also when the demand is uncertain. The online procedure reduces the impact of inaccurate forecasts and will be interesting to implement for an aggregator.

Several extensions of this work can be undertaken. First of all, our online procedure can be directly applied in the presence of other sources of stochasticity such as interactions  with market prices or local renewable production sources. The aggregator objective can also be generalized to take into account the distance to a reference load profile or to maximize consumption during renewable production peaks or could also reflect market prices.
 Also, two main theoretical questions are still open. First, the result on the PoA bound could be improved to be tighter to the numerical results, and generalized to a larger set of functions. Second, the convergence theorem for the Best Response  Dynamics (CBRD) could also be improved, as the observed convergence rate is faster than the given bound, and the convergence is numerically observed for a larger set of prices than affine functions.

\ifCLASSOPTIONcaptionsoff
  \newpage
\fi

\renewcommand{\bibfont}{\small}
\bibliographystyle{IEEEtranN}
\bibliography{../../USEFULPAPERS/Biblio_complete/shortJournalNames,../../USEFULPAPERS/Biblio_complete/biblio1,../../USEFULPAPERS/Biblio_complete/biblio2,../../USEFULPAPERS/Biblio_complete/biblioBooks}


%

\iflongversion

\appendices

\newcommand{\mum}{\underline{\mu}}
\newcommand{\mup}{\overline{\mu}}

\section{Proof of \Cref{th:uniqueness-NE}: Uniqueness of NE in $\G$}
\label{app:uniqueness-NE}
\newcommand{\hx}{\hat{\x}}
\newcommand{\hxag}{\hat{\xag}}
\newcommand{\hmp}{\hspace{-2pt}}

The proof follows the one of \cite{orda1993competitive}, extending it to the constrained case with constraints of the form \eqref{cons:minmax-power}. 

 We denote by $\lambda_n$ the Lagrange multiplier associated to \eqref{cons:total-power}, along with $\mum\nt \geq 0$ the multiplier associated to $\underline{\x}\nt \leq \x\nt$ and $\mup\nt\geq 0$ the multiplier associated to $\x\nt \leq \overline{\x}\nt$.

Note that the KKT conditions give that, at optimality:
\begin{equation} \label{eq:optim-condition-KKT}
\gamma\nt(\x\nt,\xag\ti)=\lambda_n + \mum\nt - \mup\nt \ ,
\end{equation}
where $\gamma\nt(\x\nt,\xag\ti)\eqDef c\ti(\xag\ti) +\x\nt c\ti'(\xag\ti)$ is the marginal cost of $\n$.
Let's consider $\bm{\x}$ and $\hat{\bm{\x}}$ two NE's. From \eqref{eq:optim-condition-KKT}, we  get:
\begin{align*}
& \x\nt < \overline{\x}\nt \Rightarrow \mup\nt=0 \Rightarrow  \gamma\nt(\x\nt,\xag\ti) \geq \lambda_n  \\ 
\text{ and }  & \x\nt > \underline{\x}\nt \Rightarrow \mum\nt=0  \Rightarrow  \gamma\nt(\x\nt,\xag\ti) \leq \lambda_n
\end{align*}
and the same inequalities hold for $\hat{\x}$.

We start by showing the following implications:
\begin{align}
\label{eq:implic hatx lower x} 
\left(\hat{\lambda}_n \leq \lambda_n \text{ and } \hxag\ti \geq \xag\ti \right) \ \Rightarrow \hat{\x}\nt \leq \x\nt \ , \\ 
\label{eq:implic hatx greater x}
\left(\hat{\lambda}_n \geq \lambda_n \text{ and } \hxag\ti \leq \xag\ti \right) \ \Rightarrow \hat{\x}\nt \geq \x\nt \ . 
\end{align}
We show \eqref{eq:implic hatx lower x} as \eqref{eq:implic hatx greater x} is symmetric. If $\hx\nt\hmp=\hmp\underline{\x}\nt$ or $\x\nt\hmp=\hmp\overline{\x}\nt$, then $\hx\nt \leq \x\nt$ is clear. Else, $\hx\nt >\underline{\x}\nt$ and $\x\nt<\overline{\x}\nt$ so:
\hspace{-5pt }\begin{equation} \label{eq:marginal cost ineq 2 NE}
\hmp\gamma\nt(\hx\nt,\hxag\ti) \leq \hat{\lambda}_n \hmp\leq\hmp \lambda_n \hmp\leq \hmp \gamma\nt(\x\nt,\xag\ti) \leq \hmp\gamma\nt(\x\nt,\hxag\ti)
\end{equation}
where the last inequality holds as $\gamma\nt$ is increasing in $\xag\ti$. As $c'\ti(\hxag\ti)>0$ from \Cref{ass:basic-properties-prices}, we deduce that $\x\nt  \geq \hx\nt$.

Now, let's consider $\T_1=\left\{\t \ : \ \hxag\ti > \xag\ti \right\}$ along with $\T_2=\T\setminus \T_1=\left\{ \t \ ; \ \hxag\ti \leq \xag\ti \right\}$ and $\mathcal{N}_a= \left\{n \ ; \ \hat{\lambda}_n > \lambda_n \right\}$. Suppose $\T_1 \neq \emptyset$. From constraint \eqref{cons:total-power} and from \eqref{eq:implic hatx greater x}, we have:
\begin{equation*}
\forall n \hspace{-2pt} \in \mathcal{N}_a,  \sum_{\t\in\T_1} \hspace{-5pt}\hx\nt =E_n -\hspace{-3pt}\sum_{\t\in\T_2}\hspace{-5pt}\hx\nt \leq E_n -\hspace{-3pt}\sum_{\t\in\T_2}\hspace{-5pt}\x\nt =\hspace{-3pt}  \sum_{\t\in\T_1}\hspace{-5pt}\x\nt   \ .
\end{equation*}
On the other hand, considering for $\t\in\T_1$ and $n \notin \mathcal{N}_a$, we have from \eqref{eq:implic hatx lower x} that $\hx\nt \leq \x\nt$, and thus:
\begin{equation}
\sum_{\t\in\T_1} \hxag\ti = \sum_{\t\in\T_1} \sum_{n\in \mathcal{N}_a} \hx\nt + \sum_{\t\in\T_1} \sum_{n\notin \mathcal{N}_a} \hx\nt  \leq  \sum_{\t\in\T_1} \xag\ti \ ,
\end{equation}
which is a contradiction. 
Thus $\T_1=\emptyset $ and $\forall \t, \ \hxag\ti= \xag\ti$.

With this equality, we can precise  \eqref{eq:implic hatx lower x} with:
\begin{align}
\label{eq:implic hatx lower strict x } 
\begin{split}
& \hspace{-12pt}\left[\hat{\lambda}_n\hspace{-3pt}<\hspace{-3pt}\lambda_n\hspace{-2pt}\text{ and }\hspace{-2pt}\hxag\ti\hspace{-3pt}=\hspace{-3pt}\xag\ti \right]\hspace{-2pt}\hspace{-1pt}\Longrightarrow 		\\ & \hspace{-3pt}\left[\hat{\x}\nt\hspace{-3pt}<\hspace{-3pt}\x\nt  \hspace{-2pt}\text{ or }\hspace{-2pt} \hx\nt\hspace{-3pt}=\hspace{-4pt}\x\nt\hspace{-3pt}=\hspace{-3pt}\underline{\x}\nt\hspace{-2pt}\text{ or }\hspace{-2pt}\hx\nt\hspace{-3pt}=\hspace{-3pt}\x\nt\hspace{-3pt}=\hspace{-3pt}\overline{\x}\nt \right]
\end{split}
\end{align}
and similarly for \eqref{eq:implic hatx greater x}.
Indeed, if $\hx\nt= \underline{\x}\nt$ or if $\x\nt = \overline{\x}\nt$  then the implication holds because $\x\nt \geq \lbx\nt$ and $\x\nt \leq \ubx\nt$. Else, $\hx\nt > \lbx\nt$ and $\x\nt < \ubx\nt$, and the same sequence of inequalities as in \eqref{eq:marginal cost ineq 2 NE} gives $\gamma\nt(\hx\nt,\xag\ti) <\gamma\nt(x\nt,\xag\ti)$, implying that $\hx\nt <\x\nt$.

Suppose that there exists $n$ s.t. $\hat{\lambda}_n < \lambda_n $. If only the two latter cases in \eqref{eq:implic hatx lower strict x } happen, then we have $\x\nt=\hat{\x}\nt, \ \forall \t$. Else, there is at least one $\t$ for which $\hx\nt <\x\nt$, so $E_n= \sum_\t \hx\nt < \sum_\t \x\nt=E_n$ which can not happen. Thus, $\hat{\lambda}_n=\lambda_n$ for all $n$ and \eqref{eq:implic hatx lower x} and \eqref{eq:implic hatx greater x} imply that $\x\nt=\hx\nt$ for all $n$ and $\t$.

\section{Proof of \Cref{th:bound_PoA}: PoA upper bound} \label{app:proof-PoAbound}
\newcommand{\+}{\hspace{-3pt}+\hspace{-3pt}}
\renewcommand{\=}{\hspace{-3pt}=\hspace{-3pt}}

\renewcommand{\b}{\kappa}
\newcommand{\yag}{y}
\newcommand{\y}{y}
\newcommand{\bx}{\bm{x}}
\newcommand{\by}{\bm{y}}
\newcommand{\bxag}{x} 
\renewcommand{\yag}{y}

\noindent The proof relies on the notion of \emph{local smoothness} introduced in \cite{roughgarden2015local}. The idea is to get a tighter bound than \cite{roughgarden2015local} by specifying the parameters of the  affine price functions $(c_t)_t$ and by using the upper bound on $\xag_t$ instead of looking at the worst possible cases as done in \cite{roughgarden2015local}.

Let $\b\ti\eqDef {\ca_\t}/{\cb_\t} $ so $c_\t(\bxag)=\cb_\t(\bxag+\b)$.
From \cite{roughgarden2015local}, we know that if there exist $\lambda,\mu>0$ and  a profile $\by\in\X$ satisfying for each $t\in\T$: \begin{equation}\label{eq:local-smooth-affine-condition}
\forall \bx \in \X, \ \yag\ti(\bxag\ti+\b\ti)+ \frac{\yag\ti^2}{4} \leq \lambda  \yag\ti(\yag\ti+\b\ti) + \mu \bxag\ti(\bxag\ti+\b\ti) ,
\end{equation} where $\yag\ti=\sum_n \y\nt$ and $\bxag\ti=\sum_n x\nt$ then $\G$ is locally $\lambda,\mu$-smooth for $\by$, \textit{i.e.} for any admissible profile $\bx \in \X$: $$\sum_{n=1}^N b_n(\bx)+ \nabla_n b_n(\bx)^T (\by_n-\bx_n) \leq \lambda \SC(\by) +\mu \SC(\bx).$$ In that case, it follows from \cite{roughgarden2015local} that the $\PoA$ is bounded by $\lambda/(1-\mu)$. For the remaining of the proof, we fix $t$ and omit subscript $t$ in the notations. As done in \cite{roughgarden2015local}, we introduce:
\begin{align*}
& \phi_{xy}(\mu)\eqDef \frac{y(x+\b)+ \frac{y^2}{4}-\mu x(x+\b)}{y(y+\b)} \\
  \text{ and } & \lambda^*(\mu) \eqDef \sup_{x,y\geq 0} \phi_{xy} (\mu) \ .
\end{align*}
$\lambda^*(\mu)$ is the minimum value of $\lambda>0$ such that \eqref{eq:local-smooth-affine-condition} holds with values  $(\lambda,\mu)$. Let us compute an explicit expression of $\lambda^*(\mu)$.
If $x=0$, $
\phi_{0,y}(\mu)
=\frac{y+4b}{4(y+\b)}$ and $
\dpart{\phi_{0,y}}{y} 
<0$ so $\sup_{x,y}\phi_{x,y}$ would be attained with $y=0$ and is $\phi_{0,0}=1$. Otherwise:
\begin{align*}
0=\frac{\partial \phi}{\partial x}\Leftrightarrow\frac{1}{y(y+\b)}(y-2\mu x - \mu \b)\Rightarrow x  =\frac{y-\b\mu}{2\mu} 
\end{align*}
but as $x\geq 0$, this supposes that $y \geq \mu \b$. We compute:
\begin{equation*}
\phi_{\frac{y-\b\mu}{2\mu},y} =\frac{1}{y(y+\b) 4\mu} \left( (y+\b\mu)^2 + \mu y^2 \right) \eqDef h(y) \ .
\end{equation*}
We can see that $h'$ vanishes on $\rr_+$ at $y_+\eqDef\frac{\b\mu^2 + \b\mu \sqrt{\mu^2+1-\mu}}{1-\mu}$ that gives a min of $h$ so $h$ is decreasing then increasing.
At the lower bound $y=\b\mu$, we get $\phi=\frac{\b\mu+4b}{4(\b\mu+\b)}=\frac{\mu+4}{4(\mu+1)}=\frac{1}{4} +\frac{3}{4(\mu+1)} <1 $ which is not max as $\phi_{0,0}=1$. At the upper bound $y=\ubxa$, we have $h(\ubxa)=\frac{(\ubxa+\b\mu)^2+\mu \ubxa^2}{\ubxa(\ubxa+\b)4\mu}=\lambda^*(\mu)$. Last, to compute the best bound $\inf_{\mu} \lambda^*(\mu) / (1-\mu)$, let us consider:
\begin{align*}
g(\mu) &\eqDef 4\ubxa (\ubxa+\b) \frac{\lambda^*(\mu)}{1-\mu}=\frac{(\ubxa+\b\mu)^2 +\mu \ubxa^2}{\mu(1-\mu)} \ . 
\end{align*}
If we denote $\varphi \eqDef (1+r)^2 $ and $r\eqDef {\b}/{\ubxa}$, $g(\mu)$ is minimal at $\mu^*
\eqDef ({-1+\sqrt{1+\varphi}})/{\varphi}$. We finally get our PoA  bound   as:
\begin{align*}
& \frac{\lambda^*(\mu^*)}{1-\mu^*} \= \frac{(3+2r) \+ 2 \sqrt{1+\varphi}}{4(1+r)} \=\frac{1}{2}\left(1+\sqrt{1+\frac{1}{\varphi}}+\frac{1}{2\sqrt{\varphi}} \right) \\
& = \textstyle\frac{1}{2}\left( 1+\sqrt{1+(1+r)^{-2}}+({2(1+r)})^{-1}\right) \leq 1+ \frac{3}{4(1+r)} \ .
\end{align*}
The last inequality gives a more explicit bound and is obtained from $\sqrt{a^2+b^2}\leq a + b $ valid for any $a,b \geq 0$.
%

\renewcommand{\bx}{\xx}
\renewcommand{\xag}{L}
\renewcommand{\y}{\x'}
\renewcommand{\yag}{L'}
\renewcommand{\by}{\xx'}
\noindent Next, following \cite{roughgarden2015local},  for $\xx,\by\in \X $ (admissible solutions):
\begin{align}
\nonumber
&\sum_{n} b_n(\bx) + \nabla_n b_n(\bx)^T ( \by-\bx) \\
\nonumber
&= \sum_{n} \sum_{\t\in \T} \x\nt \cdot c_\t(\xag\ti) + (\y\nt-\x\nt)\left(c_\t(\xag\ti) +\x\nt c_\t'(\xag\ti) \right) \\
\nonumber
&=\sum_\t \yag\ti \cdot c_\t(\xag\ti) + c_\t'(\xag\ti) \sum_\n\left( \y\nt \x\nt- {\x\nt}^2 \right) \\
\nonumber
&\leq \sum_\t \yag\ti \cdot c_\t(\xag\ti) \+ c_\t'(\xag\ti) \cdot \frac{ {\xag\ti}^2 }{4} = \sum_\t b_\t \left[ \yag\ti (\xag\ti + \b_\t) \+ \dfrac{{\xag\ti}^2 }{4} \right] \\
\label{eq: lambda-mu ineq hours}
& \leq \sum_\t b_\t \left[  \lambda {\yag\ti } ( \yag\ti+ \b_\t) + \mu \xag\ti (\xag\ti+\b_\t) \right] \\
\nonumber
&= \lambda \SC(\by) + \mu \SC(\bx)
\end{align}
where \eqref{eq: lambda-mu ineq hours} is valid if $(\lambda  ,\mu)$ is chosen such that: $$\forall \t \in \T ,  \  \lambda \geq  \frac{(\ubxa\ti+\b_\t \mu)^2+\mu \ubxa\ti^2}{\ubxa\ti(\ubxa\ti+\b_\t)4\mu} \eqDef \lambda^*_{\b_\t}(\mu) \ . $$
 Let us denote $\t_0 \eqDef \argmin{\t}\ \b_\t$  and choose $\mu^*\eqDef \mu_{\t_0}^*$, $\lambda^*\eqDef \lambda^*_{\b_{\t_0}}(\mu^*) $ (the optimal $(\lambda,\mu)$ for $\t_0$), then we have to check that for all $\t \in \T$, $\lambda^* \geq  \lambda^*_{\b_{\t}}(\mu^*)$. For that, it is  sufficient to show that $r \mapsto \lambda^*_r(\mu^*)$ is decreasing on $[r_{\t_0},r_\t]$, which is true if $ r_\t < -1+\sqrt{1+\frac{1-\mu^*}{{\mu^*}^2  }} \Longleftrightarrow \varphi_{r_\t} < \varphi_{r_{\t_0}}+2+\sqrt{1+\varphi_{r_{\t_0}}}$ with $\varphi_{r}=(1+r)^2$, which  
gives condition \eqref{ass-prices-for-PoAbound} stated in \Cref{th:bound_PoA}.

\section{Proof of \Cref{prop:stable game HLP}: Strong Stability of $\G$}
\label{app:strong-stability}
\renewcommand{\+}{\hspace{-3pt}+\hspace{-3pt}}
\renewcommand{\-}{\hspace{-3pt}-\hspace{-3pt}}
\renewcommand{\r}{r}
%
%
%

We denote by $G(\xx) \eqDef J_F(\xx) $ the Jacobian of operator $F$. Since functions $b_n$ are twice differentiable, condition \eqref{eq:stability} is equivalent to having the matrix $G(\xx)+G^T(\xx)$ positive definite for all $\xx \in \X$, that is, $G(\bm{\x})+G(\bm{\x})^T\succ 0$.

 As $b_n= \sum_\t b\nt$, with $b\nt(\xx\ti)\eqDef \x\nt c_\t(\xag\ti)$, is separable in $t$, we can re-index the matrix $G(\xx)$ to have a diagonal block hourly matrix $G(\bm{\x})=diag(G_1,...G_H)$ , with $G_\t(\xx\ti) \eqDef \begin{pmatrix}
 \dfrac{\partial^2 b\nt}{\partial \x\nt \partial \x\mt}(\bm{\x}\ti)
\end{pmatrix}_{n,m\in\mathcal{N}^2}$  and we get for all $t$:
\begin{equation*}
G_\t(\xx\ti)+G_\t(\xx\ti)^T=\begin{pmatrix}
 \dfrac{\partial^2 b\nt (\bm{\x}\ti)}{\partial \x\mt \partial \x\nt } + 
  \dfrac{\partial^2 b\mt (\bm{\x}\ti)}{\partial \x\nt \partial \x\mt}
\end{pmatrix}_{n,m}.
\end{equation*}
Let $\t\in \T $ and $\bm{x}\in \mathbb{R}^N  \setminus \{0\}$. Furthermore, let $\sigma(\bm{x,\x}) \eqDef\bm{x}^T\left(G_\t(\xx\ti)+G_\t^T(\xx\ti)\right)\bm{x}$. For notation simplicity, let us forget the index $t$ and the argument $(\xag)$ in functions $c_t$. We have:
%
%
\begin{align*}
\hspace{-5pt}\sigma&(\bm{x,\x})\hspace{-2pt}= \hspace{-4pt} \sum_{n=1}^N \hspace{-2pt} 2x_n^2(\x_nc'' \+ 2c') \+ 2 \hspace{-3pt} \sum_{n<m} \hspace{-3pt}x_nx_m \left( (\x_n \+ \x_m) c'' \hspace{-3pt}+\hspace{-3pt}2c'\right) =\\
& \hspace{-6pt}\sum_{n=1}^N \hspace{-3pt} 2 x_n^2 \left( \r_n\gamma \+ (1\-\r_n)a  \right)  \+ 2\hspace{-3pt}\sum_{\substack{n<m}}\hspace{-3pt} \x_n x_m \left( (\r_n\+\r_m)\gamma \+(1\-\r_n\-\r_m)a  \right)
\end{align*}
with $r_n\eqDef \x_n/\xag$, $a=2c'$ and $\gamma\eqDef 2c' +\xag c''$. Then we have:
\begin{align*}
\sigma &\= a \sum_n \hspace{-3pt}x_n^2 \+  a \Big( \sum_n \hspace{-3pt} \left(\hspace{-2pt} 1\-\r_n \left( 1\-\textstyle\frac{\gamma}{a} \right) \hspace{-3pt} \right)\hspace{-2pt}\x_n\Big)^2\hspace{-3pt} \-\textstyle\frac{\left(a -\gamma \right)^2 }{a}\hspace{-3pt}\displaystyle\sum_{n,m}\hspace{-3pt} \r_n\r_m x_nx_m 
\end{align*}
which is the sum of three quadratic form:  $q_1(\bm{x})=2a\sum x_n^2$ which has one eigenvalue $2a$ of multiplicity $N$,  $q_2(\bm{x})=  2a\left( \bm{x^T v^Tvx}\right)$ with \ $v_n \eqDef \sum_n 1-\frac{\x_n}{\xag} \left( 1-\frac{\gamma}{2a} \right)  $ of rank one which has eigenvalues $0$ of multiplicity $N-1$ and $2a||\bm{v}||_2^2$ of multiplicity $1$, and a negative form of rank one $q_3(\bm{x})=- \frac{1}{2a}\left(2a-\gamma \right)^2 \left( \sum_{n,m} \frac{\x_n}{\xag} \frac{\x_m}{\xag}\x_nx_m \right)$ which nonzero eigenvalue  is $- \frac{1}{2a}\left(2a-\gamma \right)^2 \sum_n \left(\frac{\x_n}{\xag}\right)^2$. 

We deduce that the quadratic form $q_1+q_2$ is positive definite, and that its eigenvalues are $2a$ with multiplicity $N-1$ and $2a(1+||\bm{v}||_2^2)$ with multiplicity 1.
Next we use the following result from perturbation theory:
\begin{theorem}[{\citauty{horn2012matrix}, \cite[p367]{horn2012matrix} }]

If $A,E \in \mathcal{M}_n$ are two Hermitian matrices and if $\lambda_1 \leq ... \leq \lambda_n$ are the ordered eigenvalues of $A$ and $\hat{\lambda}_1 \leq .. \leq \hat{\lambda}_n$ are the ordered values of $A+E$, and $\lambda_1^E \leq ... \leq \lambda_n^E$ are the ordered eigenvalues of $E$, then :
\begin{align*}
\forall k=1..n, & \ \lambda_1^E \leq \hat{\lambda}_k-\lambda_k \leq \lambda_n^E \\
\text{ and } &  \left| \hat{\lambda}_k-\lambda_k \right| \leq  \rho(E) = ||| E |||_2 \ .
\end{align*}
\end{theorem}
Applying this theorem with $A=q_1+q_2$ and perturbation $E=q_3$
 we get that the smallest eigenvalue $\hat{\lambda}_1$ of $ \sigma$ verifies: 
\begin{align*}
\hat{\lambda}_1 &\geq \min \left\{\text{Sp}(q_1\+q_2)\right\} \- \textstyle\frac{\left(a-\gamma \right)^2}{a} \displaystyle\sum_n r_n^2  \\
&= a \Big( 1 \- \left(1\-\textstyle\frac{\gamma}{a}\right)^2 \sum_n \r_n^2  \Big).
\end{align*}
Replacing $a$ and $\gamma$, we can get the condition of \Cref{ass:strong-stability-condition}.

%

\section{Proof of \Cref{th:cvg-BR-linear}: Convergence of CBRD}
\label{app:convergence-CBRD}
The key of the proof is that, under \Cref{ass:PricesLinear}, the game is an exact potential game \cite{monderer1996potential} with convex potential:
$$\Phi(\xx)= \sum_{\t\in\T}  \ca\ti \xag\ti+ \frac{\cb\ti}{2} ( \xag\ti^2 + \sum_n \x\nt^2 ) \ , $$
that is, for any $\xx\in\X$ and any $n$, $\nabla_n \Phi(\xx)= \nabla_n b_n(\xx)$.
Thus, the $\NE$ corresponds to the minimum of $\Phi$ and we have, for any $\xx\in\X$, $\argmin{\xx_n\in\X_n} \ b_n(\xx_n,\xx_{-n}) =\argmin{\xx_n\in\X_n} \ \Phi (\xx_n,\xx_{-n}) $. Therefore, running \Cref{algo:BRD} is equivalent to performing an alternating block coordinate minimization on $\Phi$. According to \cite[Thm.~6.1]{hong2017iteration}, we get that:
\begin{equation} \label{eq:convg-potential}
\Phi(\xx^{(k)})-\Phi(\xx^\NE) \leq \frac{2 M N^2 R^2 c_5}{k}
\end{equation}
with $M=\max_n M_n = 2 \max_t \cb\ti$ (max of Lipschitz constants of $\nabla_nb_n= \nabla_n \Phi$),
$R= \max_{\xx} \{ \norm{\xx-\xx^\NE } \ ; \Phi(x) \leq \Phi(\x^{(0)}) \}$ and $c_5= \max\{\frac{2}{MN^2R^2}-2, \ \Phi(\xx^{(1)})-\Phi(\xx^\NE), \ 2 \}$.
But $\Phi$ is also strongly convex, that is, for 	any $\xx, \xx' \in \X$:
\begin{equation} \label{eq:strong-cvx}
\Phi(\xx)-\Phi(\xx') \geq \langle \nabla \Phi(\xx') , \xx - \xx' \rangle + \frac{\stronga}{2} \norm{\xx-\xx'}^2
\end{equation}
with $a= 2\min_t \cb\ti$. Also, the minimality of $\xx^\NE$ on the convex set $\X$ implies that for any $ \xx\in \X$:
\begin{equation}\label{eq:optim-condition}
\langle \nabla \Phi(\xx^\NE), \xx - \xx^\NE  \rangle \geq 0 \ .
\end{equation}
Then from \eqref{eq:strong-cvx} and \eqref{eq:optim-condition}, we get for any $k\geq 0$:
\renewcommand{\-}{\hspace{-3pt}-\hspace{-3pt}}
\begin{align*}
\frac{a}{2}\norm{ \xx^{(k)}\-\xx^\NE}^2  & \leq \Phi(\xx^{(k)})\- \Phi(\xx^\NE) +  \langle \nabla \Phi(\xx ^\NE) , \xx^\NE \- \xx^{(k)} \rangle \\
& \leq  \Phi(\xx^{(k)})- \Phi(\xx^\NE) \ ,
\end{align*}
and from \eqref{eq:convg-potential} we get the convergence result of \Cref{th:cvg-BR-linear}.
\section{Proof of \Cref{th:SIRD-convergence}: Convergence of SIRD}
\label{app:convergence-SIRD}
\renewcommand{\-}{\hspace{-3pt}-\hspace{-3pt}}
\newcommand{\opc}{\mathrm{T}_\gamma}
\renewcommand{\bx}{\xx}
\renewcommand{\y}{\x'}
\renewcommand{\by}{\xx'}
\newcommand{\mx}{\X}

We analyze the convergence of the sequence $(\opc^k(\bm{x}))_{k}$ where $[\opc(\bx)]_n\eqDef\Pi_{\X_n}\left(\bx-\gamma \nabla_n f_n(\bx_n,\bx_{-n}) \right)$. 
First notice that the unique NE of the game $\bx^{\NE}$is the unique fixed point of $\opc$ i.e. $\bm{x}^{\NE}=\opc(\bm{x}^{\NE})$. The idea is then to prove that $\opc$ is a $\eta$-contraction, for a given norm $\Vert\cdot \Vert$ which will imply the convergence rate 
$$
\Vert T^k\big(\bm{x}^{(0)}\big)-\bm{x}^{\NE}\Vert\leq \eta^k \Vert \bm{x}^{(0)}-\bm{x}^{\NE}\Vert\ ,
$$
for any initial condition $\bm{x}^{(0)}\in \mathbb{R}^m$. 
Let $\Vert\cdot\Vert$ denote the Euclidean norm on $\mathbb{R}^d$ for any positive integer $d$. As the projection on a convex set is nonexpansive \cite[Corollary 12.20]{rockafellar2009variational}, we get for $\bx, \by \in \X$:
\begin{align*}
\label{eq:Tcontract}
\hspace{-0.5cm}&\Vert  \opc(\bx)-\opc(\by)\Vert ^2  =\sum_{\n=1}^N \Vert T_{\gamma,\n}(\bx)-T_{\gamma,\n}(\by)\Vert^2  \\
=\sum_{\n=1}^N &\Vert \Pi_{\mx_\n}(\bx_\n-\gamma \nabla_\n f_\n(\bx))- \Pi_{\mx_\n}(\by_\n-\gamma \nabla_\n f_\n(\by))\Vert^2 \nonumber \\
\hspace{-6pt}\leq\sum_{\n=1}^N & \Vert \bx_\n\-\by_\n\+\gamma (\nabla_\n f_\n(\by)\- \nabla_\n f_\n(\bx))\Vert^2 \nonumber \\
& = \sum_{\n=1}^N  \Vert \bx_\n\-\by_\n\Vert ^2  +\gamma^2  \Vert \nabla_\n f_\n(\bx)\-\nabla_\n f^h_\n(\by)\Vert^2 \\ & -2\gamma \left\langle \nabla_\n f_\n(\bx)\-\nabla_\n f_\n(\by)  ,  \bx_\n\-\by_\n) \right \rangle .  
\end{align*}

As we assume that for any $\n$, $\nabla_{\n}f_\n$ is $\lipCons_\n$-Lipschitz and let $\lipCons\eqDef\max_\n \lipCons_\n$, then we have $
\sum_{\n=1}^N \vert \nabla_\n f_\n(\bx)-\nabla_\n f_\n(\by)\vert^2\leq N \lipCons^2 \Vert \bx-\by\Vert^2$. Besides, from $\stronga$-strong stability \Cref{def:strongly-stable-game}, we get :
\begin{align*}
\Vert \opc(\bx)-\opc(\by)\Vert^2&\leq \eta \Vert \bx-\by\Vert^2\ ,
\end{align*}
with $ \eta \eqDef 1+N\lipCons^2 \gamma^2-2\gamma \alpha $. Minimizing on $\gamma>0$ gives 
 $\gamma = \frac{\alpha}{N\lipCons^2}$ and $\eta=1-\frac{\alpha^2}{N\lipCons^2}< 1$  and $\opc$ is a contraction.

\section{Proof of \Cref{th:consistency-DR-proc}: consistency of DR Procedure} \label{app:proof-DR-procedure}
\renewcommand{\xx}{\bm{x}}
\renewcommand{\x}{x}
\newcommand{\yy}{\bm{y}}
\renewcommand{\ll}{\bm{\ell}}
\renewcommand{\tend}{T}
\renewcommand{\tst}{1}
Let $t_0 \in \{\tst,\dots,\tend-1 \} $ and let us denote by $\G^{(t_0)}$ the DR-game on hours $\{t_0, \dots,\tend\}$ (considered at $t_0$ in the procedure) and $\G^{(t_0+1)}$ the DR-game on hours $\{t_0+1,\dots \tend\}$. Let  $(\xx,\bm{\lambda},\bm{\mum},\bm{\mup})$ be the  unique NE of $\G^{(t_0)}$ associated with its dual variables (as defined in \Cref{app:uniqueness-NE})  and $\yy$ the NE of $\G^{(t_0+1)}$.  As we assumed perfect forecasts of the nonflexible load, the price functions $(c_t)_t$ considered in the BR of the players do not depend on the forecasts and are the same for the games  $\G^{(t_0)}$  and $\G^{(t_0+1)}$ . We want to show that $\xx\ti=\yy\ti$ for any $t\in \{t_0+1,\dots,\tend\}$. As there is a unique NE of $\G^{(t_0+1)}$ , it is sufficient to show that $(\xx,\bm{\lambda},\bm{\mum},\bm{\mup})$ verifies the KKT conditions of each subproblem \eqref{eq:user_problem} of $\G^{(t_0+1)}$ . 
As $\xx$ is a feasible solution of $\G^{(t_0)}$ , it verifies the power bounds constraints \eqref{cons:minmax-power} as well as the total energy constraint \eqref{cons:total-power}:
\begin{equation}
\x_n^{t_0} + \sum_{t>t_0}\x\nt = E_n \ \Leftrightarrow \ \sum_{t>t_0}\x\nt = (E_n-x_n^{t_0}) \  , 
\end{equation}
where this last equality is exactly the total energy constraint \eqref{cons:total-power} for $n$ in the game $\G^{(t_0+1)}$. Also, $\xx,\bm{\lambda} ,\bm{\mum},\bm{\mup}$ verify the complementarity constraints associated to constraint  \eqref{cons:minmax-power} at each time $t \in \{t_0+1,\dots, \tend\}$. Finally, the stationarity condition gives for any $t \in \{t_0+1, \dots,\tend\},  \ \gamma\nt(\x\nt,\xag\ti)=\lambda_n + \mum\nt - \mup\nt  $ with $\gamma\nt$ the marginal cost of $n$ for time $t$. As each problem is convex, KKT conditions characterize the solution of each user's optimization problem \eqref{eq:user_problem}, and thus $\xx$ is an NE of $\G^{(t_0+1)}$, and is therefore equal to $\yy$.

\fi

%
%

%

%
%
%
%
%




\end{document}